\documentclass[12pt]{article}

\usepackage[dvips]{graphicx}
\usepackage{latexsym}
\usepackage{amssymb}

\newcommand{\CP}{\mathbb{CP}}
\newcommand{\RP}{\mathbb{RP}}
\newcommand{\R}{\mathbb{R}}
\newcommand{\T}{\mathbb{T}}

\newcommand{\D}{\mathbb{D}}

\renewcommand{\d}{\mathrm{d}}

\newcommand{\koniec}{\begin{flushright}  $\Box $ \end{flushright}}
\def\be{\begin{equation}}
\def\ee{\end{equation}}

\def\Sm{\Sigma}
\def\v{\bf v}
\def\u{\bf u}

\def\OO{\cal O}
\def\om{\omega}

\def\p{\partial}

\def\ll{\lambda}

\newtheorem{theo}{Theorem}[section]

\newtheorem{defi}[theo]{Definition}

\begin{document}
\title{Oriented straight lines and twistor correspondence}

\author{Maciej Dunajski\\
Department of Applied Mathematics and Theoretical Physics, \\
University of Cambridge,\\
Wilberforce Road,\\
Cambridge CB3 0WA, UK
}
\date{November 15, 2004} 
\maketitle
\abstract{The tangent bundle to the $n$--dimensional  sphere
is the space of oriented lines in $\R^{n+1}$. We characterise 
the smooth sections of $TS^n\rightarrow S^n$ which correspond 
to points in $\R^{n+1}$
as gradients of eigenfunctions of the Laplacian on $S^n$ with
eigenvalue $n$. The special case of $n=6$ and its connection with 
almost complex geometry is discussed.} 
\vskip5pt
\noindent
\section{Oriented lines in $\R^{n+1}$}
Oriented geodesics in $\R^{n+1}$ are straight lines. They can be parametrised
by choosing a unit vector $\u$ giving a direction, and taking the position
vector $\v$ of the point on the geodesic nearest to the chosen origin.
A pair of vectors $(\u, \v)$ corresponds to the oriented line ${\v}+t\u$, where
$t\in\R$. The space of oriented geodesics is then given by
\be
\label{twistor_sp}
{\T}=\{({\u}, {\v} )\in S^{n}\times\R^{n+1}, \;{\u.\v} =0 \}.
\ee 
For each fixed $\u$ this space restricts to a tangent plane to a unit 
$n$--sphere,
and so ${\T}$ is just the tangent bundle $TS^{n}$. We shall call
${\T}$ the twistor space.
There exists a fix-point-free map $\tau:{\T}\longrightarrow{\T}$,
such that $\tau^2=1$, obtained by reversing the orientation of each geodesic,
i.e. $\tau(\u, \v)=(-\u, \v)$.

Let $p$ be a point in $\R^{n+1}$ with a position vector ${\bf p}$. 
The oriented lines through ${p}$ are parametrised by the unit $n$--sphere in 
$T_p\R^{n+1}$,
and therefore each $p$ corresponds to a section 
$L_p:S^{n}\longrightarrow TS^{n}$ given by
\be
\label{quad_sec}
\u\longrightarrow ({\u, {\bf s(u)}), \qquad\mbox{where}\qquad  
{\bf s(u)}   ={\bf p}}-{({\bf p}.\u)\u}.
\ee
Note that these sections are preserved by $\tau$.
Each section vanishes at two points, where $\pm{\bf p}$ is 
normal to the sphere.
\subsection{Laplace sections}
The Euclidean group $E(n+1)$ acts on $\R^{n+1}$ and on $TS^{n}$, and
$E(n+1)/so(n+1)=\R^{n+1}$, so the preferred sections are orbits of $so(n+1)$.
The ($n+1)$-dimensional space of preferred sections of $\T$
corresponding to points in 
$\R^{n+1}$ can be characterised as  the eigenspace of the 
Laplacian on the n-sphere with eigenvalue $n$,   
with the vector fields being the  gradients for the eigenfunctions.
\begin{defi}
\label{laplace_sec}
The gradients of eigenfunctions of the Laplacian on $S^n$ with eigenvalue 
$n$ are called the Laplace sections of $TS^n$.
\end{defi}

\begin{theo}
\label{twistor_theorem}
There is a one-to-one correspondence between (Fig. \ref{Ntwistor_fig}) 
\begin{figure}
\caption{Twistor Correspondence}
\label{Ntwistor_fig}
\begin{center}
\includegraphics[width=8cm,height=5cm,angle=0]{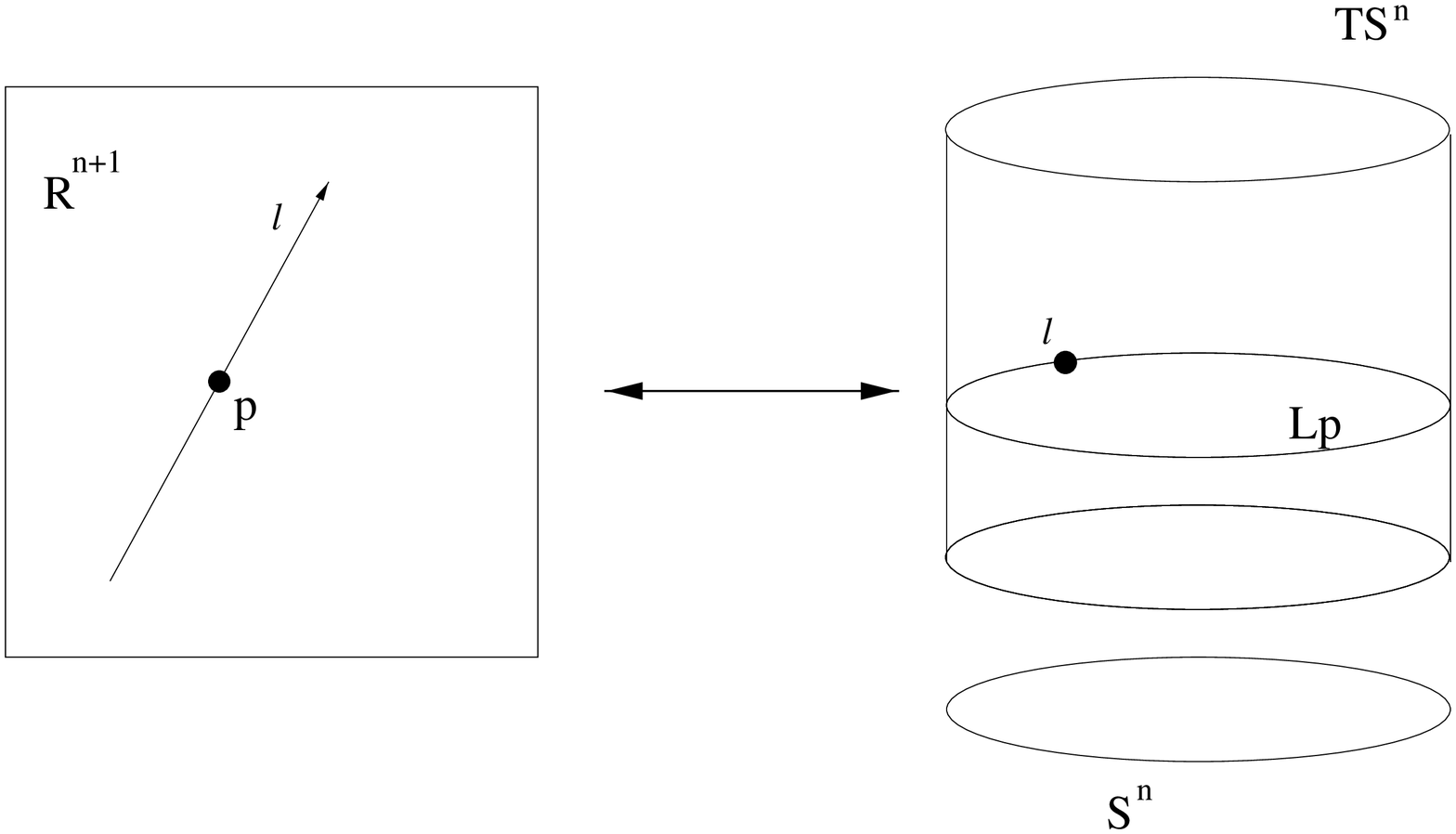}
\end{center}
\end{figure}
\begin{eqnarray*}
 \R^{n+1}&\longleftrightarrow& TS^n\\
\mbox{Points} &\longleftrightarrow& \mbox{Laplace sections} \\
\mbox{Oriented lines}&\longleftrightarrow&\mbox{Points.}
\end{eqnarray*}
\end{theo}
{\bf Proof.} To complete the proof we need to show that all Laplace
sections are of the form (\ref{quad_sec}) in some coordinates.
To see it consider a unit sphere $S^n$ isometrically immersed in
$\R^{n+1}$,
and identify a point of $S^n$ with a unit position vector $\u$. 
Let $h$ be the Riemannian metric on $S^n$ induced
by the Euclidean inner product on $\R^{n+1}$, and
let $X, Y\in T_{\u} S^{n}$. Then 
\[
{\nabla'_X}Y=\nabla_X Y+h(X, Y)\u
\]
where $\nabla'$ is the flat connection on $\R^{n+1}$, and $\nabla$ is 
the induced connection on the sphere.
If $F:\R^{n+1}\rightarrow \R$, then
\be
\label{lap_formula}
\triangle_{\R^{n+1}}(F)=-r^{-n}\frac{\p}{\p r}\Big(r^n\frac{\p F}{\p r}\Big)
+r^{-2}\triangle_{S^n}(F|_{S^n(r)}),
\ee
where $\triangle_{\R^{n+1}}=-\nabla'\cdot \nabla'$ is the Laplacian on 
$\R^{n+1}$, $\triangle_{S^n}$ is the Laplacian on the unit $n$--sphere,
and  $F|_{S^n(r)}$ is the restriction of $F$ to an $n$--sphere of radius $r$.

For any constant vector 
${\bf p}\in \R^{n+1}$ consider a function $\chi({\bf u})=\u\cdot{\bf p}$
on $S^n$. We verify that
\be
\label{gradient}
\nabla'(\chi)=-r^{-1}({\bf p}-
({\u}\cdot{\bf p}){\u}), 
\qquad  \triangle_{\R^{n+1}}(\chi)=\frac{n}{r^2}\chi.
\ee
Restricting the Laplacian to the unit sphere with $r=1$ we deduce that
\be
\label{laplacian}
\triangle_{S^n} ({\chi})=n\chi.
\ee
In particular each coordinate
function in $\R^{n+1}$ regarded  as a function on $S^n$ is an eigenfunction
of $\triangle_{S^n}$ with an eigenvalue $n$. 

The space of solutions
to (\ref{laplacian}) is $n+1$ dimensional and the 
bijection between linear functions on $\R^{n+1}$ and solutions to
(\ref{laplacian}) can be established as follows:
We have already verified that restrictions of linear functions from
$\R^{n+1}$ to $S^n$ satisfy (\ref{laplacian}). Conversely, let $\chi:S^n\longrightarrow \R$ satisfy  (\ref{laplacian}). Using the representation 
(\ref{lap_formula}) we deduce that $r\chi$ is 
a harmonic function homogeneous of degree one on $\R^{n+1}$. Let $x_i$ be local
coordinates on $\R^{n+1}$ with $|x|=r$. Therefore for each $i=1, ..., n+1,$ 
$\p (r\chi)/\p x_i$ is harmonic and homogeneous of degree $0$, and so it descends to a harmonic function on $S^n$. There are no such functions apart 
form the constants, so we deduce that $r\chi={\bf x}\cdot{\bf p}$, thus
establishing the 
bijection\footnote{Another (equivalent) characterisation of the preferred sections (\ref{quad_sec}) is a
direct consequence of 
(\ref{gradient}). 
Consider the infinitesimal generators ${\bf s}$ of non-homothetic conformal 
transformations, such that  ${\bf s}=\nabla\chi$. The equation 
${\cal L}_sh=2\chi h$ will then imply that $\chi$ satisfies 
(\ref{laplacian}).}.
\koniec
This argument can be extended to show that 
the space of homogeneous harmonic polynomials on $\R^{n+1}$ of degree $k>1$,
when restricted to $S^n$ constitute the eigenspace of $\triangle_{S^n}$
with eigenvalue $k(k+n-1)$. The multiplicity of this eigenvalue is
(consult \cite{BGM71} for details)
\[
{n+k\choose k}-{n+k-2 \choose k-2}.
\]

Let us list the properties of the Laplace sections which follow
from Theorem \ref{twistor_theorem}
\begin{itemize}
\item Laplace sections are invariant under a map  
$\tau:TS^n\longrightarrow TS^n$ given by reversing orientations
of lines in $\R^{n+1}$.
\item Each non-zero Laplace section vanishes at exactly two points
on $S^n$. Two distinct non-zero Laplace sections $L_p$ and $L_q$ 
intersect at two points in $TS^n$. These points correspond to two 
oriented lines joining $p,q \in \R^{n+1}$
\[
\u=\pm\frac{{\bf p}-{\bf q}}{|{\bf p}-{\bf q}|},\qquad 
{\bf v}=\frac{{\bf p}\cdot{\bf q}-|{\bf q}|^2}{|{\bf p}-{\bf q}|}{\bf
  p}+\frac{{\bf p}\cdot{\bf q}-|{\bf p}|^2}{|{\bf p}-{\bf q}|}
{\bf q}.
\]
Three (or more) Laplace sections generically don't meet.
\end{itemize}
To make the whole  construction independent on the choice of the 
origin in $\R^{n+1}$, we should regard the twistor space as an affine
vector bundle over $S^n$ with no preferred zero section.  

The twistor space $\T$ can also be obtained by
factoring the correspondence space $S^n\times\R^{n+1}$ 
by the action $({\u}, {\bf v})\longrightarrow  ({\u}, t{\u}+ {\bf v})$ for
$t\in\R$. This action is generated by the geodesic
flow $X$, and leads to a double fibration
$$
\begin{array}{rcccl}
&&S^n\times\R^{n+1}&&\\
&p_2\swarrow&&\searrow p_1&\\
&{\R^{n+1}}&&\T&
\end{array}
$$
given by
\[
p_2({\u}, {\bf v})={\bf v}, \qquad
p_1({\u}, {\bf v})=(\u, {\bf v}-({\bf v}\cdot {\u}){\u}).
\]

Let us look at some special cases: (here $\nabla=\p/\p\u$)
\begin{itemize}
\item For $n=1$ the unit circle $S^1$ is parametrised by $\phi\in[0,
  2\pi]$, ${\bf p}=(x_1, x_2)$, and
\[
{\bf s(u)}\cdot
\nabla=\mbox{Re}\;\Big((x_1+ix_2)\exp{(i\phi)}\frac{\d}{\d \phi}\Big).
\]
\item For $n=2$ one easily verifies
\[
{\bf s(u)}\cdot \nabla=\mbox{Re}\;\Big(((x_1+ix_2)+2\ll x_3-\ll^2(x_1-ix_2))
\frac{\d}{\d \ll}\Big),
\] 
where $\ll =(u_1+iu_2)/(1-u_3)$ is a holomorphic coordinate on $\CP^1=S^2$,
and ${\bf p}=(x_1, x_2, x_3)$.
The Laplace sections are in this case holomorphic sections of
$T\CP^1$ preserved by $\tau$. This is the original twistor
correspondence established by Hitchin \cite{H82} in his construction
of magnetic monopoles, and recently used in \cite{GK04} in a study of
generalised surfaces in ${\R^3}$. 

A much older application goes back to Whittaker \cite{W03}. We shall explain
it in a modern language of Hitchin:
Given an 
element of  $f\in H^1(T\CP^1, {\OO}(-2))$) restrict it to a Laplace 
section. The general
harmonic function on $\R^3$ is then given by
\[
V(x_1, x_2, x_3)=\oint_{\Gamma} f(\ll, (x_1+ix_2)+2\ll x_3-\ll^2(x_1-ix_2))\d \ll,
\]
where $\Gamma\subset L_p\cong\CP^1$ is a real closed contour.

A different integral transform (the X-ray transform introduced by 
John \cite{J38}) can be used to construct solutions to 
ultra-hyperbolic wave equation on the twistor space. This takes a smooth
function on $\RP^3$ (a compactification of $\R^3$) and integrates it 
over an oriented geodesic. 
The resulting function is defined on the Grassmannian 
$\mbox{Gr}_2(\R^4)$ of two-planes in $\R^4$ and satisfies the 
wave equation for a flat metric in $(++--)$ signature.
\end{itemize}
\section{Almost complex structure and $TS^6$}
The Riemannian connection $\nabla$ on $S^n$ can be used to define an
almost complex structure on $TS^n$ for any $n$. Let $T(TS^n)=V\oplus H$
be the splitting of the tangent space to $TS^n$ into vertical and horizontal 
components. Define $J_D: TS^n\longrightarrow TS^n$ by
\[
J_D(X_H)=X_V, \qquad J_D(X_V)=-X_H,
\]
where $X_V$ and  $X_H$ are the vertical and horizontal parts of a vector on
$TS^n$. This structure was studied by Dombrowski \cite{D62} who showed that
the torsion of $J_D$ does not vanish unless both the torsion and the
curvature of $\nabla$ are zero. This almost complex structure has
nothing to do with the Laplace sections defined in Def. 
\ref{laplace_sec}. 
From now on
we shall restrict to the case $n=6$ where another (inequivalent)
almost complex structure can be defined on $\T$.  The basic 
facts about the cross products on $\R^7$ will be
recalled, and used to show that the Laplace sections are almost complex. 
\subsection{Cross product in $\R^7$ and the group $G_2$}
Let $(x_1, ..., x_7)$ be coordinates on $\R^7$, and let
$\d x_{ijk}$ be a shorthand notation for $\d x_i\wedge \d x_j\wedge \d x_k$.
Following Bryant \cite{B87} we define the exceptional group $G_2$ as
\[
G_2=\{\rho\in GL(7, \R)| \rho^*(\phi)=\phi\},
\]
where 
\[
\phi=\d x_{123}+ \d x_1\wedge(\d x_{45}+\d x_{67})+\d x_2\wedge(\d x_{46}-\d x_{57})
-\d x_3\wedge(\d x_{47}+\d x_{56}).
\]
It is a compact, connected, and simply connected Lie group of dimension $14$.
It also preserves the Euclidean metric 
$g=\d x_1^2+...+\d x_7^2$, the orientation $\d x_{1234567}$, and the four-form
\[
*\phi= \d x_{4567}+ \d x_{23}\wedge(\d x_{45}+ \d x_{67}) 
-\d x_{13}\wedge(\d x_{46} -\d x_{57}) 
-\d x_{12}\wedge(\d x_{47}+\d x_{56}). 
\]
The group $G_2$ acts transitively on a unit sphere $S^6\subset \R^7$ with
a stabiliser $SU(3)$.

A cross product $\times:\R^7\times\R^7\longrightarrow \R^7$ can be defined by
\[
g(X\times Y, Z)=\phi(X, Y, Z).
\]
This cross product has the same properties as the one induced by
the octonion multiplication, which leads to a more standard definition of 
$G_2$ as the group of automorphisms of the octonions. The induced  cross 
product satisfies the identities analogous to those in three-dimensions
\be
\label{identities7}
g(X\times Y, X\times Y)=g(X, X)g(Y, Y)-g(X, Y)^2,
\qquad
X\times(X\times Y)=g(X, Y)X-g(X, X)Y.
\ee
\subsection{Pseudoholomorphic sections of $TS^6$}
Consider a curve $\gamma(s, t)$ of oriented 
lines in $\R^7$ parametrised by $s\in\R$, and given by
\[
{\bf\gamma}(s, t)={\v} (s)+t{\u} (s).
\]
A ${\u}$-orthogonal projection of 
tangent vector $t\dot{\u}+\dot{\v}$ gives rise to a normal Jacobi field
\be
\label{split_jacobi}
V=(\dot{\v}-(\dot{\v}. {\u}){\u}+t\dot{\u})|_{s=0}
=(\dot{\u}, \dot{\v}-(\dot{\v}. {\u}){\u}),
\qquad\mbox{where}\qquad
\dot{\;}=\frac{\p}{\p s}.
\ee
All vectors tangent to a space of oriented geodesics are of this form.

Let us define a map $\widetilde{J}:T{\T}\longrightarrow T{\T}$ by
\be
\label{def_structure}
V\longrightarrow \widetilde{J}(V)={\u}\times V,\qquad\mbox{where}\qquad 
V\in T_{({\u, \v})}\T.
\ee
From the properties (\ref{identities7}) of 
cross-product $\times$ in $\R^7$ it follows that 
$\widetilde{J}$ is an almost complex structure. Indeed,
\[
\widetilde{J}^2(V)={\u}\times({\u}\times V)=({\u .}V){\u}-({\u .\u})V=-V.
\]
Note that $\tau(\widetilde{J})=-\widetilde{J}$.

This almost complex structure is related to a standard almost 
complex structure $J$ on $S^6$ defined by
$J(\v)=\u\times \v$.
To see this
consider the restriction of the Euclidean scalar product from $\R^7$
to $S^6$. This gives the unique nearly K\"ahler metric $h$ on $S^6$ compatible
with $J$ \cite{FI55} in a sense that 
\[
h(X, Y)=h(JX, JY), \qquad  \nabla_{X}J(X)=0, \qquad \forall X, Y\in TS^6, 
\]
where $\nabla$ is the Levi--Civita connection of $h$. 
Let
\[
T(TS^6)=V\oplus H
\]
be the splitting of the tangent space to $TS^6$ into vertical and horizontal 
components with respect to $\nabla$. The almost complex
structure on $TS^6$ defined   
by taking the standard almost complex structure $J$ on each factor $H$ and $V$
coincides with the almost--complex structure 
(\ref{def_structure}), because the splitting (\ref{split_jacobi})
coincides with the  splitting $T(TS^6)$ induced by $\nabla$ (which is
a projection of splitting given by restricting $\nabla'$ to a tangent space).
In particular $\widetilde{J}$ is not integrable, since $J$ isn't.

Let $\rho:S^6\rightarrow S^6$ be an element of $G_2$, and let 
${\v}\in T_{\u} S^6$. Then
\[
\rho_*(J({\v}))=\rho({\u})\times\rho_*({\v})=\widetilde{J}(\rho_*({\v}))\in T_{\rho({\u})}S^6.
\]
and we deduce that the Laplace 
sections $L_p$ of ${\T}\longrightarrow S^6$
which correspond to points in $\R^7$ are $G_2$--invariant in a sense that
\[
\rho_*(L_{p}({\u}))=L_{\rho({\bf p})}\rho(\u).
\]

Now we want to argue that the Laplace 
sections are also almost complex in the sense that
\[
\widetilde{J}\circ (L_p)_*=(L_p)_*\circ J.
\]
This follows directly from the
geometrical construction because $\u$ is a unit normal to a sphere of 
geodesics $L_p$ through $p$, and the cross product preserves the almost
complex structure on $S^6$ (the almost complex structure on the space
of lines  is a rotation in $\R^7$ through $90$ degrees about the
direction of the line which preserves the tangent spaces of $L_p$). 

It can also be seen by applying $\widetilde{J}$ to (\ref{split_jacobi}) and
performing a direct calculation. This leads to an overdetermined
system of equations for $L:S^6\rightarrow TS^6, L({\u})=(u^j, L^j(u))$
\[
\Big(\phi_{ljm}u^j\Sm_{pk}+\phi_{kjp}u^j\Sm_{ml}\Big)\frac{\p L^m}{\p
u^p}=0,
\]
where $\Sm_{ij}=\delta_{ij}-u_iu_j$. These equations are satisfied by
the Laplace sections.
\section{Other twistor correspondences}
In this final section we shall mention two other generalisations
of the Hitchin correspondence. 
The first one (due to Study \cite{S03} for $n=2$)
is more than hundred years old. The second one (due to Murray
\cite{M85}) gives a way of solving the Laplace equation.

{\bf Study's correspondence.}
The correspondence between oriented lines in $\R^{n+1}$ and points in $TS^n$ 
can be re-expressed in terms of the dual numbers of the form
\[
a+\tau b
\]
where $a, b \in \R$, and $\tau^2=0$. 
Let $\D$ denote the space of the dual numbers. Any oriented line in $\R^{n+1}$
can be represented by a vector in $\D^{n+1}$ 
\[
{\bf A}=\u +\tau \v
\]
which is of unit length with respect to an Euclidean norm in $\D^{n+1}$ induced from 
$\R^{n+1}$. This gives an analogue of Study's result \cite{S03}:
There is a one to one correspondence between oriented lines in $\R^{n+1}$ and
points on the dual unit sphere in $\D^{n+1}$. Comparing 
this with (\ref{twistor_sp}), we see that the dual unit sphere
in $\D^{n+1}$ is equivalent to $TS^n$ with an additional structure (that of 
dual numbers) selected on the fibres. 

Let $\theta$ and $\rho$ be the angle and the distance between two oriented 
lines represented by ${\bf A}$ and ${\bf B}$. 
Define a dual angle by 
\[
\Theta=\theta+\tau \rho.
\]
Using a formal definition
\[
\cos{\Theta}=1-\frac{1}{2!}\Theta^2+\frac{1}{4!}\Theta^4 +...= 
\cos{\theta}-\tau\sin{\theta},
\]
one can verify an attractive looking formula 
\[
{\bf A}\cdot{\bf B}=\cos{\Theta},
\]
and deduce that group of Euclidean motions in $\R^{n+1}$ is equivalent to
$O(n+1, \D)$. 

 {\bf Murray's correspondence.} Let 
$[z_0, z_1, ..., z_n]$ be homogeneous coordinates on $\CP^n$, and
let $f=z_0^2+z_1^2+...+z_n^2$ define a section of ${\OO}(2)\longrightarrow \CP^n$.
This section vanishes on a hyper-quadric 
\[
X=\{f=0, [z]\in \CP^n\}\subset\CP^n.
\]
Murray \cite{M85} defines a  twistor space $Z$ 
to be a restriction of the total space of
${\OO}(1)\longrightarrow \CP^n$ to $X$. This leads to a double fibration
\[\begin{array}{rcccl}
&&X\times\R^{n+1}&&\\
&m_2\swarrow&&\searrow m_1&\\
&{\R^{n+1}}&&Z.&
\end{array}
\]
The canonical bundle of $K_X$ of $X$ in $Z$ is ${\OO}(-n+1)$.
\begin{theo}[Murray \cite{M85}] Let $\triangle_{\R^{n+1}}$ be the Laplacian
on $\R^{n+1}$. There exists an isomorphism
\[
T: H^{n-1}(Z, K_X)\longrightarrow Ker\; (\triangle_{\R^{n+1}})
\]
given by
\[
T(\om)(z)=\int_{X_z}\om,
\]
where $(\om)$ is a $K_X$-valued $(0, n)$ form on $Z$ pulled back to 
$X\times\R^{n+1}$.
\end{theo}
The twistor spaces $\T$ and  $Z$ have the same dimensions, but
the connection between Theorem \ref{twistor_theorem} 
and the Murray correspondence is not clear.
\section*{Acknowledgements}
I thank Michael Eastwood, Nigel Hitchin and Simon Salamon
for useful discussions, and  Marc Lachi\`eze-Rey for
pointing out some errors in an earlier version of this paper.

\end{document}